\documentclass[12pt]{amsart}
\usepackage{amsthm}
\usepackage{amssymb}
\usepackage{amsmath}

\newtheorem{lemma}{Lemma}[section]
\newtheorem{theorem}[lemma]{Theorem}

\newtheorem{proposition}{Proposition}
\newtheorem*{corollary*}{Corollary}

\numberwithin{equation}{section}
\newcommand{\st}[1]{\ensuremath{^{\scriptstyle \textrm{#1}}}}

\gdef\myletter{}

\let\savetheequation\theequation
\def\theequation{\savetheequation\myletter}

\def    \Cinf   {C^\infty}

\newcommand{\CC}{{\mathbb C}}
\newcommand{\LL}{{\mathbb L}}
\newcommand{\RR}{{\mathbb R}}

\newcommand{\fg}{{\mathfrak g}}

\renewcommand{\div}{\mbox{div}}
\renewcommand{\Im}{\mbox{Im}}

\newcommand{\eff}{\mbox{eff}}
\newcommand{\Log}{\mathop{\mathrm{Log}}}
 \newcommand{\red}{\mbox{red}}
\renewcommand{\Re}{\mbox{Re}}   
\newcommand{\Span}{\mbox{span}}
\newcommand{\stable}{\mbox{stable}}
\newcommand{\vol}{\mbox{vol}}

\renewcommand{\O}{{\mathcal O}}

\def \bar{\overline}
\def \hat{\widehat}

\begin{document}

\title[Potential functions and actions of tori]{Potential functions
and actions of tori on K\"ahler manifolds}  

\author{D. Burns and V. Guillemin}\thanks{Supported in part by NSF
grants DMS-0104047 (DB) and DMS-0104116 (VG)}

\address{M.I.T., Cambridge, MA, {\em and} Univ. of Michigan, Ann
Arbor, MI}  
\email{dburns@umich.edu}

\address{M.I.T., Cambridge, MA}
\email{vwg@math.mit.edu}

\date{September 24, 2002}



\maketitle


\section{Introduction}
\label{sec:1}

Let $M$ be a K\"ahler manifold equipped with a free Hamiltonian
action of the standard $n$-torus, $T$ with moment map, $\Phi : M
\to \RR^n$.  For $\lambda \in \RR^n$ the symplectic quotient
\begin{displaymath}
  M_{\lambda} = \Phi^{-1} (\lambda) /T
\end{displaymath}
inherits from $M$ a K\"ahler structure; and in the first part of
this paper we will describe what the K\"ahler form and Ricci form
look like locally on coordinate patches in $M_{\lambda}$.  Then
in the second part of this paper we will discuss some global
implications of these results.  This will include

\begin{enumerate}
\item 
A K\"ahlerian proof of the Duistermaat-Heckman theorem.

\item 
A formula, due to Biquard and Gauduchon, for the K\"ahler
potential on a symplectic quotient.

\item 
  A convexity theorem of Atiyah for moment images of
  $T_{\CC}$-orbits.

\item 
  A formula in terms of moment data for the K\"ahler metric on a
  toric variety.

\item 
  A formula for the K\"ahler form on the symplectic quotient of a
  K\"ahler--Einstein manifold.
\end{enumerate}

\noindent
A few comments about each of these items:

\begin{enumerate}
\item 
The usual proof of the Duistermaat-Heckman theorem is global and
topological in nature.  Our proof in the K\"ahler case consists of
showing that locally on a coordinate patch in $M_{\lambda}$ two
canonically defined K\"ahler forms are identical.

\item 
The formula of Biquard--Gauduchon was used by Calderbank, David
and Gauduchon to give an elegant economical proof of the theorem
alluded to in item~4.  (This theorem is due to one of us, and a
much less elegant proof of it was presented in \cite{Gu}.)  In
this paper we use a primitive local form of the
Biquard--Gauduchon theorem (see \S~4, formula~(4.5)) to give an
even simpler derivation of this formula.

\item 
  Atiyah's proof of this convexity theorem for $T_{\CC}$-orbits
  involves a careful analysis of the global geometry of the closures
  of these orbits.  The proof we present here is quite different in
  spirit.  In the spirit of this article it is based on convexity
  properties of strictly pluri-subharmonic potentials. Unfortunately,
  we must assume that the K\"ahler form defines an integral cohomology
  class, an hypothesis Atiyah avoids. Although standard approximation
  techniques will eliminate this hypothesis in some cases, we are
  unable to derive Atiyah's theorem in full generality. We note that
  K\"ahler potentials have already been used to prove convexity
  theorems for more general moment maps by P. Heinzner and
  A. Huckleberry \cite{HH}.

\item 
It is well-known that the symplectic quotient by a circle action of a
K\"ahler--Einstein manifold is not in general K\"ahler--Einstein.
(See, for instance, \cite{Be} and \cite{Fu}\S~7.3
\footnote{We take this opportunity to remark that the examples of
non-existence of extremal metrics on ruled surfaces cited in \cite{Fu}
as ``Burns's examples'' are due to P. DeBartolomeis and one of us \cite{BD}}.)
 As was pointed out to us by Michael Duff, this is analogous to
 what happens, vis a vis dimensional reduction by circle actions,
 in general relativity.  By Kaluza--Klein, the quotient by a
 circle action of a $4+1$-dimensional Einstein manifold has, in
 addition to its downstairs metric, a two-form (magnetic field)
 and an electrostatic potential; and the field equations satisfied
 by these data are not just the field equations for gravity
 alone, but for gravity \emph{coupled} with electromagnetism.  In
 other words the field equations on the quotient space reflect
 the fact that the geometry of the quotient space
is more
 complicated than the geometry of the original space.

We will show that this is true for symplectic quotients as well.
For $\lambda \in \RR^n$, the level set, $\Phi^{-1} (\lambda)$, is
a principal $T$-bundle over $M_{\lambda}$, and the metric on $M$
gives one a connection on this bundle and a set of curvature
forms, $\mu_i$.  In addition one gets from the metric a potential
function on $M_{\lambda}$, the effective potential:  $V_{\eff}$,
whose value at $p \in M_{\lambda}$ is the volume of the fiber
above $p$ in $\Phi^{-1} (\lambda)$.\footnote{For circle actions,
  this coincides with the usual ``effective potential'' for
  reduced systems in classical mechanics.  (See
  \cite{AM},\S~4.5.)}  We will show that if the K\"ahler form,
$\omega$, and the Ricci form, $\mu$, satisfy on $M$ the
K\"ahler--Einstein equation, $\mu = \kappa \omega$, then the reduced
K\"ahler form, $\omega_{\lambda}$, and reduced Ricci form,
$\mu_{\lambda}$, satisfy
\begin{equation}
  \label{eq:1.1}
  \mu_{\lambda} - 2\sqrt{-1} \partial \, \bar{\partial} \Log
  V_{\eff} + \sum c_i \mu_i 
  = \kappa (\omega_{\lambda} + \sum \lambda_i \mu_i)
\end{equation}
the $c_i$'s being constants which don't depend on $\lambda$.
This identity can be thought of as a K\"ahler version of the
Einstein equation for a metric coupled with an electromagnetic
field.

The presence of the $c_i$'s in (\ref{eq:1.1}) is due,
incidentally, to the fact that the moment map is only
well-defined up to an additive constant.  We will show in
\S~\ref{sec:11} that, for one very natural normalization of this
constant, the $c_i$'s
 are equal to zero.  The proof of this will depend on a fact
 about K\"ahler--Einstein manifolds which is of independent
 interest.  Let $Z_i$ , $i=1,\ldots ,n$, be the complex vector
 fields generating the action of $T_{\CC}$  on $M$.  Then the
 components, $\phi_i$, of the moment map satisfy the divergence
 identites
 \begin{equation}
   \label{eq:1.2}
   \kappa \phi_i = \div Z_i +c_i \, .
 \end{equation}
We remark that equations like \ref{eq:1.1} have been exploited earlier
by C. LeBrun and coworkers for the construction of K\"ahler-Einstein
and constant scalar curvature surfaces: see especially \cite{L} and {PP}.

\end{enumerate}

\section{The Darboux theorem for $T$-invariant K\"ahler forms}
\label{sec:2}

Let $\CC^* = \CC - \{ 0 \}$ and let $T_{\CC} = (\CC^*)^n$ be the
complexification of the group, $T=(S^1)^n$.  Let $U$ be an open
convex subset of $\CC^m$ and let
\begin{equation}
  \label{eq:2.1}
  M = T_{\CC} \times U = (\CC^*)^n \times U
\end{equation}
be the trivial $T_{\CC}$-bundle over $U$.  We will prove in this
section an equivariant Darboux theorem for K\"ahler forms on $M$.
Let $\omega$ be a $T$-invariant K\"ahler form which is Hamiltonian
with respect to the action of $T$; and let $z_1,\ldots ,z_n$ and
$w_1 , \ldots , w_m$ be the coordinate functions on $(\CC^*)^n$
and $U$.  We will prove that
\begin{equation}
  \label{eq:2.2}
  \omega = \sqrt{-1} \partial \, \bar{\partial} \rho
\end{equation}
where
\begin{equation}
  \label{eq:2.3}
  \rho = \rho (t_1,\ldots ,t_n \, , \, w_1, \ldots ,w_m), 
  \quad t_i = |z_i |^2 \, .
\end{equation}
Moreover, we will show that if $\rho_1$ and $\rho_2$ are two such
functions
\begin{equation}
  \label{eq:2.4}
  \rho_2 - \rho_1
 = \sum \lambda_i \Log t_i + \Re f
\end{equation}
with $\lambda_i \in \RR$ and $f \in \O (U)$.

\begin{proof}
  Since the action of $T$ is Hamiltonian the orbits of $T$ are
  isotropic.  However, for $p_0 \in M$, the orbit, $T \cdot
  p_0$ is a deformation retract of $M$.  Therefore, since the
  restriction of $\omega$ to this orbit is zero, $\omega $ is
  exact; i.e.,~$\omega = d\alpha$.  Let $\alpha = \alpha_{0,1} +
  \alpha_{1,0}$. Then, $\omega = \partial \, \alpha_{0,1} +
  \bar{\partial} \alpha_{1,0}$; and $\bar{\partial}
  \alpha_{0,1}=\partial \, \alpha_{1,0}=0$. Since $M$ is a domain of
  holomorphy, $\alpha_{0,1} = \bar{\partial}h$ and hence
  \begin{displaymath}
    \omega = \partial \, \bar{\partial} (h-\bar{h})
       = \sqrt{-1} \partial \, \bar{\partial} \rho
  \end{displaymath}
where $\sqrt{-1} \rho = h-\bar{h}$.  By averaging with respect to
$T$ we can assume that $\rho$ is $T$-invariant and hence of the
form (\ref{eq:2.3}).  To prove the assertion (\ref{eq:2.4}) we
make use of the following elementary fact.

\begin{lemma}
  \label{lem:2.1}
  If $u \in \Cinf (\CC^*)$ is a  radially symmetric harmonic
  function it is of the form $c_1 \Log r + c_2$, $c_i \in \CC$.
\end{lemma}

\begin{proof}
By hypothesis
\begin{displaymath}
  \Delta u = \frac{d^2}{dr^2} u+ \frac{1}{r} \, \frac{du}{dr} =0 \,
  ,
\end{displaymath}
and all solutions of this ODE are of the form above.

\end{proof}

We will now prove (\ref{eq:2.4}) by induction on $n$.  The
difference, $g=\rho_2 - \rho_1$, is harmonic in $z_1$, so it has
to be of the form
\begin{displaymath}
  g_1(z',w) \Log |z_1|^2 + g_2 (z',w), 
  \quad z'=(z_2 , \ldots , z_n) \, .
\end{displaymath}
However,
\begin{displaymath}
  0= \frac{\partial^2}{\partial \, z_1 \partial \, \bar{z}_i} g
  = \frac{1}{z_1} \, \frac{\partial}{\partial \, \bar{z}_i}
  g_1 (z',w) \, ,
\end{displaymath}
so the derivative of $g_1$ with respect to $\bar{z}_i$ is zero
and, since $g_1$ is real-valued, so is the derivative with
respect to $z_i$.  Similarly the derivatives of $g_1$ with
respect to $w_i$ and $\bar{w}_i$ are zero; so $g_1$ is constant.
Thus
\begin{displaymath}
  g=\lambda_1 \Log |z_1|^2 + g_2 (z',w).
\end{displaymath}
Now apply induction. 

\end{proof}

 The result we have just proved can be
restated:

\begin{theorem}
  \label{th:2.2}

If $g$ is a $T$-invariant real-valued pluri-harmonic function on
$M$ then $g$ is of the form
\begin{displaymath}
  \sum \lambda_i \Log t_i + \Re f
\end{displaymath}
with $\lambda_i \in \RR$ and $f \in \O (U)$.
\end{theorem}

\section{The moment map}
\label{sec:3}

We will compute below the moment map associated with the action
of $T$ on $M$.  Let $\rho $ be the potential function
(\ref{eq:2.3}).  Letting $\alpha = -i \partial \, \rho$, one gets
from (\ref{eq:2.1}), where $M$ is as in \S2.
\begin{displaymath}
  \omega = d \alpha
\end{displaymath}
and hence with 
\begin{displaymath}
  \frac{\partial}{\partial \, \theta_j} = i \left( z_j
    \frac{\partial}{\partial \, z_j} -\bar{z}_j
    \frac{\partial}{\partial \, \bar{z}_j} \right)
\end{displaymath}
one has
\begin{displaymath}
  \iota \left( \frac{\partial}{\partial \, \theta_j} \right) \omega 
  = - d  \left( \iota \left( \frac{\partial}{\partial \,\theta_j} \right) 
    \alpha \right)
\end{displaymath}
so we get for the $j$\st{th} component of the moment map
\begin{equation}
  \label{eq:3.1}
  \phi_j = \iota \left( \frac{\partial}{\partial \,\theta_i} \right)
     \alpha = z_j \frac{\partial \, \rho}{\partial \, z_j}=
     t_j \frac{\partial \, \rho}{\partial \, t_j} \, .
\end{equation}
Unlike the symplectic form, $\omega$, (\ref{eq:3.1}) depends
upon the choice of $\rho$.  If we replace $\rho$ by the potential
function
\begin{equation}
  \label{eq:3.2}
  \rho_1 = \rho + \sum \lambda_j \Log t_i
  + \Re f \, \quad f \in \O (U)
\end{equation}
and set
\begin{displaymath}
  \alpha_1 =-i\partial \, \rho_1 = -i \partial \, \rho
    -i \sum \lambda_j \frac{dz_j}{z_j}-
    \frac{i}{2} \partial \, f
\end{displaymath}
we get for the $j$\st{th} component of the moment map
\begin{displaymath}
  \iota \left( \frac{\partial \, }{\partial \, \theta_j} \right) \alpha_1
    = \iota \left( \frac{\partial}{\partial \,\theta_j}\right)
    \alpha + \lambda_j \, , 
\end{displaymath}
i.e.,~the moment map associated with $\rho_1$ is
\begin{equation}
  \label{eq:3.3}
  \Phi_1 = \Phi + \lambda \, .
\end{equation}
Thus if $\Phi$ and $\Phi_1$ are the same, $\lambda =0$ and hence
by (\ref{eq:3.2}) $\rho_1 = \rho + \Re f$, $f \in \O (U)$.  Thus
to summarize, we proved

\begin{theorem}
  \label{th:3.1}

Let $\rho_1$ and $\rho$ be strictly pluri-subharmonic functions
of the form (\ref{eq:2.3}).  If the symplectic forms and moment
maps associated with $\rho$ and $\rho_1$ are the same then
\begin{equation}
  \label{eq:3.4}
  \rho_1 = \rho + \Re f , \quad f \in \O (U) \, .
\end{equation}

\end{theorem}

Coming back to the formula (\ref{eq:3.1}) note that if we make
the change of variables $t_i = e^{s_i}$ and set
\begin{equation}
  \label{eq:3.5}
  \rho (t_1,\ldots ,t_n ,w) = \rho(e^{s_1}, \ldots , e^{s_n},w)
  = : f (s_1 , \ldots s_n ,w)
\end{equation}
the moment map can be written in the simpler form
\begin{equation}
  \label{eq:3.6}
  \phi_j (s,w) = \frac{\partial}{\partial \, s_j} f(s,w) \, .
\end{equation}

\section{The reduced potential function}
\label{sec:4}

As above let $\rho$ be the potential function (\ref{eq:2.3}) and
let
\begin{displaymath}
  \alpha = -i \partial \, \rho =-i\sum
  \frac{\partial \, \rho}{\partial \, z_j} \, dz_j 
  + \partial_w \rho \, .
\end{displaymath}
By (\ref{eq:3.1}), $z_j \partial \, \rho /\partial \, z_j$ is equal to
$\phi_j$; hence
\begin{equation}
\label{eq:4.1}
  \alpha =-i (\sum \phi_j \, d \, \Log z_j + \partial_w \rho) \, .
\end{equation}
Let us denote by $\iota_\lambda$ the inclusion of the level set
$\phi^{-1}(\lambda)$ and by $\pi$ the projection of $T_{\mathbb{C}}
\times U$ onto $U$. The composite map $\pi \circ \iota_\lambda$ is
$T$-invariant and gives us an identification locally of $U$ with the
reduced space
$$M_\lambda = \phi^{-1}(\lambda)/T,$$
and via this identification, $w_1,\ldots,w_n$ become local holomorphic
coordinates on $M_\lambda$.

{}From the identity \ref{eq:4.1}, we get
\begin{equation}
  \label{eq:4.2}
  \iota^*_{\lambda} \alpha =- i \iota^*_{\lambda}
  \partial_w \rho + dh_{\lambda}
\end{equation}
where $h_{\lambda} =- i \sum \lambda_j \Log \iota^*_{\lambda}
z_j$.  Setting $\iota^*_{\lambda} t_j = t_j (w)$, one can
interchange $\iota^*_{\lambda}$ and $\partial_w$:
\begin{displaymath}
  \iota^*_{\lambda} \frac{\partial}{\partial w} \rho (t,w)
  = \frac{\partial}{\partial \, w}\rho (t(w),w)
  -\sum \frac{\partial \, \rho}{\partial t_j} \, 
  \frac{\partial t_j}{\partial \, w},
\end{displaymath}
and, by (\ref{eq:3.1}), rewrite this in the form
\begin{equation}
  \label{eq:4.3}
  \iota^*_{\lambda} \frac{\partial}{\partial \, w} \rho (t,w)
  = \frac{\partial}{\partial \, w}
    (\iota^*_{\lambda} \rho - \sum \lambda_j 
    \Log \iota^*_{\lambda} t_j) \, .
\end{equation}
Hence (\ref{eq:4.2}) becomes 
\begin{equation}
  \label{eq:4.4}
  \iota^*_{\lambda} \alpha = -i \partial_w \iota^*_{\lambda}
   (\rho - \sum \lambda_j \Log t_j)
  + dh_{\lambda} \, ,
\end{equation}
so the reduced symplectic form $d\iota^*_{\lambda}\alpha$ on
$M_\lambda$ is equal to
\begin{eqnarray}
  \label{eq:4.5}
  \omega_{\lambda} &=& i \partial_w\,  \bar{\partial}_w \rho_{\lambda}\\
\noindent{\nonumber{where}}\\
  \label{eq:4.6}
\rho_{\lambda} &=& \iota^*_{\lambda} (\rho - \sum \lambda_j
\Log t_i) \, .
\end{eqnarray}

\section{The reduced K\"ahler metric on $M_{\lambda}$}
\label{sec:5}

By equation (\ref{eq:4.5})
\begin{eqnarray}
  \label{eq:5.1}
  \omega_{\lambda} =i \sum \psi_{\alpha ,\bar{\beta}}
    \, dw_{\alpha} \wedge \, d\bar{w}_{\beta}\\
\noindent{\nonumber{where}}\\
\label{eq:5.2}
   \psi_{\alpha,\bar{\beta}}
       =\frac{\partial^2 \, \rho_{\lambda}}{\partial \, w_{\alpha}
         \partial \, \bar{w}_{\beta}} \, .
\end{eqnarray}
To compute these derivatives we will, as in (\ref{eq:3.5}), set
$t_i = e^{s_i}$ and let
\begin{displaymath}
  f(s_1,\ldots ,s_n,w) = \rho (t_1 ,\ldots ,t_n,w) \, .
\end{displaymath}
In these coordinates the $\lambda$-level set of the moment map is
defined by the equations
\begin{displaymath}
  \frac{\partial \, f}{\partial \, s_i} = \lambda_i
\end{displaymath}
so if we let
\begin{equation}
  \label{eq:5.3}
  g(s,w) = f(s,w)-\sum \lambda_i s_i
\end{equation}
then, by (\ref{eq:4.5}), $\rho_{\lambda}$ is the restriction of
$g$ to the set
\begin{equation}
  \label{eq:5.4}
  \frac{\partial \, g}{\partial \, s_i}=0 \, .
\end{equation}
Letting $s=s(w)$ on this set, we obtain from (\ref{eq:5.2}) and
(\ref{eq:4.5}):
\begin{displaymath}
  \psi_{\alpha ,\bar{\beta}} = \frac{\partial}{\partial \,
    w_{\alpha}} \, \frac{\partial}{\partial \, \bar{w}_{\beta}}
  (g(s(w),w)) \, .
\end{displaymath}
However, by (\ref{eq:5.4})
\begin{displaymath}
  \frac{\partial}{\partial \, w_{\alpha}}
   (g(s(w),w)) = \frac{\partial \, g}{\partial \, w_{\alpha}}
   (s(w),w)
\end{displaymath}
and hence
\begin{equation}
  \label{eq:5.5}
  \psi_{\alpha ,\bar{\beta}} =\sum \left(
    \iota^*_{\lambda} \frac{\partial^2 \, g}{\partial \, s_i \partial \, w_{\alpha}}
  \right) \frac{\partial \, s_i}{\partial \, \bar{w}_{\beta}} 
  + \iota^*_{\lambda}
  \frac{\partial^2 \, g}{\partial \, w_{\alpha} \partial \, \bar{w}_{\beta}}
  \, .
\end{equation}
Moreover, by (\ref{eq:5.4})
\begin{displaymath}
  \frac{\partial}{\partial \, w_{\alpha}} \left(
    \frac{\partial \, g}{\partial \, s_i} (s(w),w) \right)
  =0= \sum_j \frac{\partial^2 \, g}{\partial \, s_i \partial \, s_j}
  (s(w),w) + \frac{\partial^2 \, g}{\partial \, s_i \partial \, w_{\alpha}}
  (s(w),w) \, .
\end{displaymath}
Hence
\begin{equation}
  \label{eq:5.6}
  \iota^*_{\lambda} \frac{\partial^2 \, g}{\partial \, s_i \partial \,
    w_{\alpha}} = -\sum \left(  \iota^*_{\lambda} 
    \frac{\partial^2 \, g}{\partial \, s_i \partial \,  s_j} \right)
  \frac{\partial \, s_j}{\partial \, w_{\alpha}}
\end{equation}
and if we substitute this into (\ref{eq:5.5}) we get the
following slightly more symmetric formula for $\psi_{\alpha
  ,\bar{\beta}}$:
\begin{equation}
  \label{eq:5.7}
  \psi_{\alpha ,\bar{\beta}} = \iota^*_{\lambda}
  \frac{\partial^2 \, g}{\partial \, w_{\alpha}\partial \, \bar{w}_{\beta}}
  -\sum \left( \iota^*_{\lambda}
    \frac{\partial^2 \, g}{\partial \, s_i \partial \, s_j} \right)
  \frac{\partial \, s_i}{\partial \, w_{\alpha}} \, 
  \frac{\partial \, s_j}{\partial \, \bar{w}_{\beta}} \, .
\end{equation}

\section{The Ricci potential}
\label{sec:6}

By definition the Ricci potential of the K\"ahler form, $i
\partial \, \bar{\partial} \rho$, is minus the log determinant of the
matrix
\begin{gather}
  \left[
    \begin{array}{cc}
      \frac{\partial^2 \, \rho}{\partial \, z_i \partial \, \bar{z}_j} &
  \frac{\partial^2 \, \rho}{\partial \, z_i \partial \, \bar{w}_{\beta}}\\[4ex]
    \frac{\partial^2 \, \rho}{\partial \, w_{\alpha}\partial \,  \bar{z}_j}     
& \frac{\partial^2 \, \rho}{\partial \, w_{\alpha} \partial \, \bar{w}_{\beta}}
    \end{array}      
\right] \, .
\tag{$*$}
\end{gather}
But with $s_i = \Log |z_i|^2$
\begin{eqnarray}
\label{eq:6.1}
  \frac{\partial^2 \, \rho}{\partial \, z_i \partial \, \bar{z}_j}
  &=& \frac{\partial^2 \,}{\partial \, z_i \partial \, \bar{z_j}}
  (\rho -\sum \lambda_{\kappa} \Log |z_{\kappa}|^2)\\[2ex]
\nonumber
  &=& \frac{\partial^2}{\partial \, z_i \partial \, \bar{z}_j}
     g(s,w) = \frac{1}{z_i \bar{z}_j} \, 
     \frac{\partial^2 \, g}{\partial \, s_i \partial \, s_j} \, .
\end{eqnarray}
Similarly
\begin{eqnarray}
  \label{eq:6.2}  
\frac{\partial^2 \, \rho}{\partial \, z_i \partial \, \bar{w}_{\beta}}
    &=& \frac{1}{z_i} \, 
    \frac{\partial^2 \, g}{\partial \, s_i \partial \, \bar{w}_{\beta}}\\[2ex]
\label{eq:6.3}  
   \frac{\partial^2 \, \rho}{\partial \, w_{\alpha} \partial \, \bar{z}_j}
    &=& \frac{1}{\bar{z}_j} \, 
       \frac{\partial^2 \, g}{\partial \, w_{\alpha}\partial \, s_j}\\
\noalign{\nonumber{and}}\\
\label{eq:6.4}
   \frac{\partial^2 \, \rho}{\partial \, w_{\alpha} \partial \,
     \bar{w}_{\beta}}
   &=& \frac{\partial^2 \, g}{\partial \, w_{\alpha}\partial \, \bar{w}_{\beta}}
\end{eqnarray}
so the determinant of the matrix above is equal to the
determinant of the matrix
\begin{displaymath}
  \left[
    \begin{array}{cc}
      \frac{\partial^2 \, g}{\partial \, s_i \partial \, s_j} &
   \frac{\partial^2 \, g}{\partial \, s_i \partial \, \bar{w}_{\beta}}\\[4ex]
      \frac{\partial^2 \, g}{\partial \, w_{\alpha}\partial \, s_j}
   & \frac{\partial^2 \, g}{\partial \, w_{\alpha} \partial \, \bar{w}_{\beta}}
    \end{array}      
\right]
\end{displaymath}
times a factor of
\begin{equation}
  \label{eq:6.5}
  \frac{1}{|z_1|^2 \cdots |z_n|^2} = e^{-(s_1+\cdots +s_n)} \, .
\end{equation}
Now consider the restriction of this matrix to the set,
$\Phi^{-1} (\lambda)$.  If we multiply the $j$\st{th} column of
this matrix by $\partial \, s_j /\partial \, \bar{w}_{\beta}$
and subtract it from the $\beta$\st{th} column we get by
(\ref{eq:5.5}) and (\ref{eq:5.6}) the matrix
\begin{displaymath}
  \left[
    \begin{array}{cc}
      \frac{\partial^2 \, g}{\partial \, s_i \partial \, s_j} &
      0 \\[2ex]
      \cdots
      & \psi_{\alpha ,\bar{\beta}}
    \end{array}      
\right] \, .
\end{displaymath}
Hence the determinant of the matrix ($*$) is
\begin{displaymath}
  e^{-(s_1 + \cdots + s_n)} \det \left( 
  \frac{\partial^2 \, g}{\partial \, s_i \partial \, s_j}\right)
\det (\psi_{\alpha ,\bar{\beta}})
\end{displaymath}
or
\begin{displaymath}
  \det \left( 
    \frac{\partial^2 \, \rho}{\partial \, z_i \partial \, \bar{z}_j}
    \right) \det (\psi_{\alpha , \bar{\beta}}) \, .
\end{displaymath}
Thus we've shown:

\begin{theorem}
  \label{th:6.1}
Let $h$ be the Ricci potential of the K\"ahler form, $\omega$.
Then the Ricci potential, $h_{\lambda}$, of the K\"ahler form
$\omega_{\lambda}$ is
\begin{equation}
  \label{eq:6.6}
  h_{\lambda} = \iota^*_{\lambda}
     \left( h + \Log \det \left(
         \frac{\partial^2 \, \rho}{\partial \, z_i \partial \, \bar{z}_j}
       \right)\right) \, .
\end{equation}

\end{theorem}

The second summand on the right has a simple geometric
interpretation.  Let $\nu_i = d \Log z_i$.  At a point $p \in M$
the restriction of $\omega$ to the orbit, $T_{\CC}\cdot p$, is
equal to
\begin{equation}
  \label{eq:6.7}
  \sqrt{-1} \sum z_i \bar{z}_j
  \frac{\partial \,^2 \rho}{\partial \, z_i \partial \, \bar{z}_j}
  \nu_i \wedge \bar{\nu}_j
\end{equation}
and the K\"ahlerian volume form on this orbit is

\begin{equation}
\label{eq:6.8}
    \sqrt{-1} f \nu_1 \wedge
    \bar{\nu}_1 \wedge \cdots \wedge \nu_n \wedge \bar{\nu}_n
\end{equation}
where
\begin{equation}
  \label{eq:6.9}
  f=2^n |z_1|^2 \ldots |z_n|^2 \det \left(
    \frac{\partial \,^2 \rho}{\partial \, z_i \partial \, \bar{z}_j}
    \right) \, .
\end{equation}
But by (\ref{eq:6.7}) the map
\begin{displaymath}
  J_p : \Re \,  \nu_i \to \Im \, \nu_i \, , \,
     i=1,\ldots ,n
\end{displaymath}
maps the space
\begin{equation}
  \label{eq:6.10}
  \Span \{ \Re \, \nu_i , \, i=1,\ldots ,n \}
\end{equation}
at $p$ isometrically onto the space
\begin{equation}
  \label{eq:6.11}
  \Span \{ \Im \, \nu_i , \, i=1,\ldots ,n \} \, .
\end{equation}
Moreover, (\ref{eq:6.7}) also implies that these spaces are
Riemannian orthocomplements of each other in the cotangent space
at $p$ to $T_{\CC} \cdot p$.  Hence if $g \Re\, \nu_1 \wedge \cdots
\wedge \Re \, \nu_n$ is the Riemannian volume form on the space
(\ref{eq:6.10}), $g \Im \, \nu_1 \wedge \cdots \wedge \Im
\, \nu_n$ is
the Riemannian volume form on the space (\ref{eq:6.11}); and the
wedge product of these two forms is the Riemannian volume form
on the cotangent space at $p$ to $T_{\CC} \cdot p$, which is just
the volume form (\ref{eq:6.8}).  Thus $g^2=f$.  Moreover, the
space (\ref{eq:6.11}) is just the cotangent space to $T \cdot p$
at $p$ and $\Im \, \nu_i = d\theta_i$.  Hence the volume form on the
$T$ orbit through $p$ is $\sqrt{f} \, d\theta_1 \wedge \cdots
\wedge d\theta_n$ and the volume of this orbit is $(2\pi)^n \sqrt{f(p)}$.
Thus $\iota^*_{\lambda} \sqrt{f}$ is the effective
potential. $V_{\eff}$, that we defined in \S~\ref{sec:1}; and, up
to an additive constant,
\begin{equation}
  \label{eq:6.12}
  h_{\lambda} =\iota^*_{\lambda} (h + \sum \Log |z_i|^2)
  + 2\Log V_{\eff} \, .
\end{equation}

\section{GIT quotients}
\label{sec:7}

In the second half of this paper we will describe some global
applications of the results of \S\S~1--6.  Let $M$ now be a
$d$-dimensional compact complex manifold, and $\tau : T_{\CC}
\times M \to M$ a holomorphic action of $T_{\CC}$ on $M$.  Let
$\omega \in \Omega^{1,1}(M)$ be a K\"ahler form with respect to
which the restriction of $\tau$ to $T$ is Hamiltonian and let
$\Phi :M \to \RR$ be the moment map, normalized so that
\begin{displaymath}
  \int \phi^i \omega^d =0 \, .
\end{displaymath}
We will assume that $0$ is a regular value of $\Phi$ and that $T$
acts freely on the level set $\Phi^{-1}(0)$.  Then, for $\lambda$
close to $0$, $T$ acts freely on $\Phi^{-1}(\lambda)$.  Hence the
reduced space
\begin{displaymath}
  M_{\lambda} = \Phi^{-1} (\lambda)/T
\end{displaymath}
is well-defined and is a compact K\"ahler manifold.  We will
denote by $\omega_{\lambda}$ its K\"ahler form.

We recall (see \cite{MFK}) that the set of \emph{stable points}
in $M$ is the set
\begin{displaymath}
  M_{\stable} = T_{\CC} \cdot \Phi^{-1}(0) \, .
\end{displaymath}
This set is an open dense subset of $M$ whose complement is a
subvariety of dimension $\leq d-1$.  Moreover, $T_{\CC}$ acts
freely and properly on $M_{\stable}$, so the quotient space
\begin{displaymath}
  M_{\red} = M_{\stable}/T_{\CC}
\end{displaymath}
is a compact complex manifold of dimension $\ell = d-n$,
and, by definition, $M_{\red}$  is the \emph{geometric invariant
  theory} (GIT) \emph{quotient of } $M$ \emph{by} $T_{\CC}$.  The
projection
\begin{equation}
  \label{eq:7.1}
  \pi : M_{\stable} \to M_{\red}
\end{equation}
makes $M_{\stable}$ into a principal $T_{\CC}$ bundle over
$M_{\red}$.  Moreover, since $\Phi^{-1}(\lambda )$ is compact it
is contained in $M_{\stable}$ for $\lambda$ close to $0$; hence
the map
\begin{displaymath}
  \Phi^{-1}(\lambda) \hookrightarrow M_{\stable} \overset{\pi}{\longrightarrow}
  M_{\red}
\end{displaymath}
induces a biholomorphic map, $M_{\lambda} \to M_{\red}$; so, as
complex manifolds, the $M_{\lambda}$'s are identical with
$M_{\red}$ and we can think of the $\omega_{\lambda}$'s as a
family of K\"ahler forms on $M_{\red}$.  In the first of our
applications we will describe how the $\omega_{\lambda}$'s vary
as one varies the parameter, $\lambda$.  In all the applications
below, by the way, except for the application discussed in
\S~\ref{sec:10}, we will only be concerned with the open subset,
$M_{\stable}$, of $M$ and the fibration (\ref{eq:7.1}).
Moreover, for the most part, it won't be important that
$M_{\red}$ be compact; in other words, in most of the
applications below, we can assume that $M=M_{\stable}$ and allow
$M$ to be a principal $T_{\CC}$-bundle over a (not necessarily
compact) base manifold, $M_{\red}$.

\section{The K\"ahler version of Duistermaat--Heckman}
\label{sec:8}

Let $U$ be a convex coordinate patch in $M_{\red}$ and let $M_U =
\pi^{-1} (U)$.  We will assume that the bundle (\ref{eq:7.1}) is
trivial over $U$ and hence
\begin{equation}
  \label{eq:8.1}
  M_U = (\CC^*)^n \times U \, .
\end{equation}
Let $\rho$ be a potential function on $M_U$ of the form
(\ref{eq:2.3}) such that $\omega = \sqrt{-1} \partial \,
\bar{\partial} \, \rho$ and such that the moment map (\ref{eq:3.1})
coincides with the moment map defined in the previous section.
By (\ref{eq:4.6}) the reduced symplectic form on $U$ is equal to
\begin{equation}
  \label{eq:8.2}
  \sqrt{-1}\,  \partial \, \bar{\partial} \iota^*_{\lambda} \rho
  - \sum \sqrt{-1} \, \lambda_j \partial \, \bar{\partial}
  \Log \iota^*_{\lambda} |z_j|^2 \, .
\end{equation}
(Here and throughout this section, we identify $T$-invariant functions
on $\phi^{-1}(\lambda)$ with functions on $U$, i.e., as functions of
the $w$ coordinates, where we can express, {\em along}
$\phi^{-1}(\lambda)$, $|z_i|^2 = |z_i|^2(w)$.)
We claim that the first term in (\ref{eq:8.2}) is an
\emph{intrinsically} defined $2$-form on $U$ not depending on the
trivialization (\ref{eq:8.1}).  In fact if $\rho_1$ is another
potential function with the same properties as $\rho$ then by
Theorem~\ref{th:3.1}, $\rho_1=\rho + \Re f$, for $f \in \O (U)$
and hence
\begin{equation}
  \label{eq:8.3}
  \sqrt{-1}\,  \partial \, \bar{\partial} \iota^*_{\lambda} \rho_1
  = \sqrt{-1} \, \partial \, \bar{\partial}\iota^*_{\lambda}\rho \, .
\end{equation}
The same argument shows that for different values, $\lambda_1$
and $\lambda_2$, of $\lambda$
\begin{equation}
  \label{eq:8.4}
  \sqrt{-1}\,  \bar{\partial} \iota^*_{\lambda_1} \rho 
  - \sqrt{-1} \, \bar{\partial} \iota^*_{\lambda_2}\rho
\end{equation}
is an intrinsically defined one-form not depending on $\rho$; for
if we replace $\rho$ by $\rho + \Re f$, this difference is
unchanged.  Thus there exists a globaly defined two-form,
$\mu_{\lambda}$, whose restriction to $U$ satisfies
\begin{equation}
  \label{eq:8.5}
  \mu_{\lambda} =\sqrt{-1}\,  \partial \, \bar{\partial}
     \iota^*_{\lambda} \rho
\end{equation}
and a globally defined one-form, $\nu_{1,2}$, whose restriction
to $U$ satisfies
\begin{equation}
  \label{eq:8.6}
  \nu_{1,2} = \sqrt{-1} \, \bar{\partial} \iota^*_{\lambda_1}
  \rho -\sqrt{-1}\, \bar{\partial}\iota^*_{\lambda_2}\rho
\end{equation}
such that
\begin{equation}
  \label{eq:8.7}
  \mu_{\lambda_1} - \mu_{\lambda_2} = d\nu_{1,2} \, .
\end{equation}

Next consider the terms in the second summand of (\ref{eq:8.2}).
The projection of the product, $T_{\CC}=(\CC^*)^n$ onto its
$i$\st{th} factor defines a character
\begin{displaymath}
  \chi_i :T_{\CC} \to \CC^* \, .
\end{displaymath}
Let $\LL_i$ be the line bundle over $M_{\red}$ associated with
this character:  A section of this line bundle is by definition a
map
\begin{displaymath}
  s:M_{\stable} \to \CC
\end{displaymath}
satisfying
\begin{equation}
  \label{eq:8.8}
  s(ap) = \chi_i (a)^{-1}s(p) \, ;
\end{equation}
so in particular the functions 
\begin{equation}
  \label{eq:8.9}
  s_i : M_U \to \CC^* \, , \, 
  (z,w) \to z^{-1}_i
\end{equation}
are non-vanishing holomorphic sections of $\LL_i$ over $U$.

Given a map $s:M_{\stable} \to \CC$ satisfying (\ref{eq:8.8}) the
restriction of $|s|^2$ to $\Phi^{-1}(\lambda)$ is, by
(\ref{eq:8.8}) constant on the fibers of the $T$-bundle:
\begin{displaymath}
  \Phi^{-1}(\lambda) \to M_{\red} \, ,
\end{displaymath}
so one can define a hermitian metric on $\LL_i$ by setting
\begin{equation}
  \label{eq:8.10}
  \langle s,s \rangle_{\lambda} = \iota^*_{\lambda}
  |s|^2 \, ,
\end{equation}
and in terms of this metric the $i$\st{th} summand of
(\ref{eq:8.2}) is equal to
\begin{equation}
  \label{eq:8.11}
  - \sqrt{-1} \partial \,\bar{\partial}
  \Log \langle s_i,s_i \rangle_{\lambda} = \mu_{i,\lambda}
\end{equation}
which is just the curvature form $\mu_{i,\lambda}$ of the bundle $\LL_i$.  This
form depends on $\lambda$, however if $\lambda_1$ and $\lambda_2$
are distinct values of $\lambda$, the inner products, $\langle \cdot ,
\cdot \rangle_{\lambda_1}$  and $\langle \cdot ,\cdot
\rangle_{\lambda_2}$ differ at each point, $w \in M_{\red}$, by a
positive multiplicative constant, $\kappa_i (w)$, and
\begin{equation}
  \label{eq:8.12}
  \mu_{i,\lambda_1}-\mu_{i,\lambda_2}
  = \sqrt{-1} \, d (\partial \, \Log \kappa_i(w)) \, ;
\end{equation}
hence $\mu_{i,\lambda_1}$ and $\mu_{i,\lambda_2}$ are canonically
cohomologous.  Hence, to summarize, the reduced symplectic form,
$\omega_{\lambda}$, on $M_{\red}$ satisfies the
Duistermaat -- Heckman equation:
\begin{equation}
  \label{eq:8.13}
  \omega_{\lambda} =\mu + \sum \lambda_i \mu_i
\end{equation}
where $\mu =\mu_{\lambda}$ is defined by (\ref{eq:8.5}) and
$\mu_{i}=\mu_{i,\lambda}$ is the curvature form of the line
bundle $\LL_i$.  These forms are \emph{canonically} defined
$2$-forms on $M_{\red}$.  Moreover, though they depend on
$\lambda$, they are, for different choices of $\lambda$,
\emph{canonically} cohomologous by (\ref{eq:8.7}) and
(\ref{eq:8.12}).

\section{The Biquard--Gauduchon formula}
\label{sec:9}

The potential function in formula (\ref{eq:8.2}):
\begin{equation}
  \label{eq:9.1}
  \rho_{\lambda} = \iota^*_{\lambda} \rho
  - \sum \lambda_i \iota^*_{\lambda} \Log  |z_i|^2
\end{equation}
is not intrinsically defined in terms of $\rho$ since the
functions, $z_i$, depend upon the choice of the trivialization
(\ref{eq:8.1}).  However, there is a simple way to make it
intrinsically defined.  The function, $\rho_{\lambda}$, lives on $U$, so
we can pull it back by the fibration (\ref{eq:7.1}) to get a
function $\pi^*\rho_{\lambda}$ on $M_U$.  We claim that the
sum,
\begin{equation}
  \label{eq:9.2}
  \hat{K}_{\lambda} = \pi^* \rho_{\lambda} +
  \sum \lambda_i \Log |z_i|^2
\end{equation}
is independent of how we trivialize the bundles, (\ref{eq:7.1})
over $U$.  If one changes this trivialization by replacing $z_i$
by $a_i (w) z_i$ then
\begin{displaymath}
  \Log |a_i (w) z_i |^2 -\pi^* \Log \iota^*_{\lambda}
  |a_i(w)z_i|^2 = \Log |z_i|^2 -\pi^* \Log \iota^*_{\lambda}
  |z_i|^2 \, ;
\end{displaymath}
so the right-hand side of (\ref{eq:9.2}) is unchanged.  Also note
that (\ref{eq:8.3}) can be written in the more instrinsic form
\begin{equation}
  \label{eq:9.3}
  \pi^* \omega_{\lambda} = \sqrt{-1} \partial \,\bar{\partial}
  \hat{K}_{\lambda \, .}
\end{equation}

In other words, if one has a K\"ahler potential, $\rho$, on $M$
satisfying $\omega = \sqrt{-1} \partial \, \bar{\partial} \rho$,
one gets an \emph{intrinsic} K\"ahler potential, $\hat{K}$, on $M$
satisfying (\ref{eq:9.3}).  Thus, even though there may not be a
globally defined K\"ahler potential for the reduced symplectic
form, $\omega_{\lambda}$, there is for its pull-back, $\pi^*
\omega_{\lambda}$, to $M$.  (See \cite{BG}, pages~291--293.)

\section{The Atiyah convexity theorem}
\label{sec:10}

Our proof of this theorem makes use of some elementary facts
about strictly convex functions:  Let $V$ be a real vector space
and $F:V\to \RR$ a smooth function.  We recall that $F$ is
\emph{strictly convex} if the Hessian $d^2F_p$, is positive
definite for all $p \in V$, and is \emph{stable} if $F(x)$ tends
to infinity as $x$ tends to infinity in $V$.  Given $\ell \in
V^*$ let $F_{\ell}(x)=F(x)-\ell (x)$.  The \emph{stability set}
of a strictly convex function, $F$, is the set, $S_F$, of all
$\ell \in V^*$ for which $F_{\ell}$ is stable.  We will need
below the following properties of this set.

\begin{proposition}
  \label{prop:1}
  $S_F$ is an open convex set.
\end{proposition}

\begin{proof}
  Let $\ell = s\ell_1 + (1-s) \ell_2$, $0 \leq s \leq 1$.  Then
  if $F_{\ell_1}(x)$ and $F_{\ell_2}(x)$ tend to $+\infty$ as $x$
  tends to infinity so does $F_{\ell}(x)$. To see that $S_F$ is open,
we may, without loss of generality assume that $\ell = 0$ and that
$x=0$ is a strict global minimum of $F$. Then we have 
$$F(x) > \min_{\{|x|=1\}} |F(x)| + c \cdot |x|, \, {\mbox{for}} |x| >
1,$$
where $c = \min_{\{|x|=1\}} \frac{\partial F}{\partial r}$, and $c >
0$ ebcause $F$ is strictly convex and $F(0)=0$. Then for all $\ell$
such that $|\ell| < c, \ell \in S_F$.

\end{proof}

\begin{proposition}
  \label{prop:2}
If $F_1$ and $F_2$ are strictly convex and $F_1-F_2$ is bounded, 
$S_{F_1} = S_{F_2}$.
\end{proposition}

\begin{proof}
$(F_1)_{\ell}(x)$ tends to $+ \infty$ as $x$ tends to infinity if
and only if $(F_2)_{\ell}(x)$ tends to $+\infty$ as $x$ tends to infinity.
\end{proof}

Let $F$ be a function of the form
\begin{equation}
  \label{eq:10.1}
  F(x) = \sum^N_{i=1} c_i e^{\alpha_i(x)}
\end{equation}
where $c_i>0$ and $\alpha_i \in V^*$.

\begin{proposition}
  \label{prop:3}
$F$ is stable if and only if zero is contained in the interior of
the convex hull of $\{ \alpha_i , i=1,\ldots ,N \}$.
\end{proposition}

\begin{proof}
If zero is contained in the exterior or boundary of this set
there exists a $\xi \in V$ such that the $\alpha$'s lie in the
half-space
\begin{displaymath}
  \{ \mu \in V^* \, , \, \mu (\xi) \leq 0 \} \, ,
\end{displaymath}
and hence $F(x)$ is bounded along the ray $x=t\xi$, $t\geq 0$.
\end{proof}

\begin{corollary*}
The stability set of the function $\Log F$ is the interior of the
convex hull of $\{ \alpha_i , \, i=1,\ldots ,N \}$.
\end{corollary*}

\begin{proof}
$\Log F-\ell = \Log Fe^{-\ell}$ where
\begin{displaymath}
  Fe^{-\ell} = \sum c_i e^{\alpha_i -\ell}\, ,
\end{displaymath}
and this is stable if and only if $\ell$ is contained in the
interior of the convex hull of the $\alpha$'s.
\end{proof}

The Legendre transform of a function, $F$, is the map
\begin{displaymath}
  L_F :V \to V^* \, , \, p \to dF_p \, .
\end{displaymath}

\begin{proposition}
  \label{prop:4}
If $F$ is strictly convex, $L_F$ maps $V$ bijectively onto $S_F$.
\end{proposition}

\begin{proof}
Let $\ell$ be in $S_F$.  Since $F_{\ell}$ is stable $F_{\ell}(x)$
tends to $+\infty$ as $x$ tends to infinity, hence $F_{\ell}$ has
at least one minimum point, $p$.  Moreover, since $F$ is strictly
convex, every critical point of $F_{\ell}$ is a non-degenerate
minimum point.  So by the Morse index theorem, $p$ is the
\emph{unique} critical point of $F_{\ell}$; and at $p$, $dF_p=\ell$.
\end{proof}

We will now prove that the Atiyah theorem is true in one very
simple special case.  Let $\tau_0$ be a \emph{linear} action of
$T_{\CC}$ on $\CC^{N+1}$ with weights, $\beta_i$, $i=1,\ldots
,N+1$ and let $\tau$ be the induced action of $T_{\CC}$ on $\CC
P^N$.  Let $p$ be the point, $[a_0 , \cdots , a_n]$, $a_0 =1$.
On the complement of the hyperplane, $z_0=0$, the canonical
K\"ahler form, $\omega_{FS}$, is equal to $\sqrt{-1} \partial \,
\bar{\partial} \rho$ where
\begin{displaymath}
  \rho = \Log \left(
\frac{|z_0|^2 + \cdots + |z_n|^2}{|z_0|^2}
\right)
\end{displaymath}
and the restriction of $\rho$ to the $T_{\CC}$ orbit through $p$
is the function
\begin{displaymath}
  \Log (c_1 t^{\alpha_1} + \cdots + c_N t^{\alpha_N}+1)
\end{displaymath}
where $\alpha_i = \beta_i - \beta_0$ and $c_i =
|a_i|^2$.  Making the substitution $t_i=e^{s_i}$, $i=1, \ldots ,n$
this becomes
\begin{equation}
  \label{eq:10.2}
  F(s) = \Log (c_1 e^{\alpha_1 (s)} + \cdots 
      + c_N e^{\alpha_N (s)}+1) \, .
\end{equation}
Now let $\Phi$ be the moment map associated with the action of
$T$ on $\CC P^N$.  By (\ref{eq:3.6}) the restriction of $\Phi$ to
$T_{\CC} \cdot p$ is just $L_F$; and, by the corollary to
proposition~\ref{prop:3}, its image is the convex hull of the set
consisting of zero
and those $\alpha_i$'s for which $c_i \neq 0$.  In particular
this proves

\begin{theorem}
  \label{th:10.1}

The moment image of $T_{\CC} \cdot p$ is a convex polytope.
\end{theorem}

We will now show that this special case of the convexity theorem
implies the convexity theorem in much greater generality:  First
of all it is clear that if $\tau_0$ is a linear action of
$T_{\CC}$ on $\CC^{N+1}$ and $\tau$ is the induced action of
$T_{\CC}$ on $\CC P^N$ then, for every non-singular projective
variety, $M \subseteq \CC P^N$, which is invariant under $\tau$,
the theorem is true for $T_{\CC}$ orbits lying in $M$.  Now
suppose that $M$ is a compact K\"ahler manifold whose K\"ahler
form, $\omega$, is an integral form.  Then $[\omega]$ is the
Chern class of a complex line bundle, $\LL$; and, by a theorem of
Kodaira,
some high power, $\LL^k$, of $\LL$ admits enough global
holomorphic sections to imbed $M$ in $\CC P^N$.  Moreover, the
action of $T_{\CC}$ on $M$ lifts to an action of $T_{\CC}$ on
$\LL$ and hence lifts to a \emph{linear} action of $T_{\CC}$ on
the space of holomorphic sections of $\LL^k$.  Thus the
Kodaira imbedding intertwines the action of $T_{\CC}$ on $M$ with
a linear action of $T_{\CC}$ on $\CC P^N$.  Let $\gamma$ be this
imbedding.  The pull-back by $\gamma$ of $\omega_{FS}$ is not the
K\"ahler form on $M$; however by Kodaira's theorem $[\omega]=\gamma^*
[\omega_{FS}]$; hence there exists a $T$-invariant function,
$H:M \to \RR$ satisfying
\begin{equation}
  \label{eq:10.3}
  \omega = \gamma^* \omega_{FS} + \sqrt{-1} \partial \, 
  \bar{\partial}H \, .
\end{equation}
Now let $q$ be a point in $M$ and let $p=\gamma (q)$.  Then by
(\ref{eq:10.3}) the restriction of $\omega$ to the orbit, $T_{\CC}
\cdot q$ is equal to $\sqrt{-1}\partial \, \bar{\partial}F_1$
where 
\begin{equation}
  \label{eq:10.4}
  F_1 (s_1 , \ldots , s_n) = H (s_1 , \ldots , s_n) +
     F(s_1 , \ldots , s_n) 
\end{equation}
$F$ being the function (\ref{eq:10.2}).  However, $H$ is bounded
on $T_{\CC} \cdot q$ so by proposition~\ref{prop:2} $L_F$ and
$L_{F_1}$ have the same image.  Thus we can conclude:

\begin{theorem}
  \label{th:10.2}
The moment image of $T_{\CC} \cdot q$  is identical with
the moment image of $T_{\CC} \cdot p$  (and hence, in
particular, by theorem~\ref{th:10.1}, is a convex polytope).
\end{theorem}

\section{The symplectic quotient of a K\"ahler--Einstein manifold}
\label{sec:11}

Suppose that the K\"ahler form on $M$ is K\"ahler--Einstein,
i.e.,~that the Ricci form and the K\"ahler form satisfy $\mu =
\kappa \omega$, $\kappa \in \RR$.  We will show that the reduced
K\"ahler form, $\omega_{\lambda}$, and reduced Ricci form,
$\mu_{\lambda}$, satisfy (\ref{eq:1.1}).  It is enough to prove this
identity on open pseudoconvex sets, $U$, of $M_{\red}$ over which
$M_{\stable}$ is the trivial bundle:  $(\CC^*)^n \times U$.  Let
$\rho$ be a potential function on $U$ of the form (\ref{eq:2.3})
such that $\omega = \sqrt{-1} \partial \, \bar{\partial} \rho$
and such that the moment map (\ref{eq:3.1}) is identical with the
given moment map.  Then if $h$ is the Ricci potential associated
with $\rho$, $\kappa \rho -h$ is pluri-harmonic so by
theorem~\ref{th:2.2}
\begin{equation}
  \label{eq:11.1}
  \kappa \rho = h+ \sum a_i \Log |z_i|^2 + \Re f \, 
\end{equation}
$f \in \O (U)$.  Thus with $c_i = (a_i-1)$ the expression
\begin{equation}
  \label{eq:11.2}
  \kappa (\iota^*_{\lambda} (\rho -\sum \lambda_i \Log |z_i|^2)
  + \sum \lambda_i \Log \iota^*_{\lambda} |z_i|^2)
\end{equation}
is equal to the expression
\begin{equation}
  \label{eq:11.3}
  \iota^*_{\lambda} (h+\sum \Log |z_i|^2 ) + 2\Log
  V_{\eff} + \Re f
\end{equation}
plus the expression
\begin{equation}
  \label{eq:11.4}
  \sum c_i \Log \iota^*_{\lambda} |z_i|^2 - 2\Log V_{\eff}
\end{equation}
Applying $\partial \, \bar{\partial}$ to these three
expressions and using (\ref{eq:6.12}) one gets (\ref{eq:1.1}). 

We will next show that the $c_i$'s are given by the divergence
equations (\ref{eq:1.2}).  Differentiating the identity
(\ref{eq:11.1}) with respect to the vector field, $Z_i=z_i
\partial /\partial z_i$, we get
\begin{equation}
  \label{eq:11.5}
  \kappa Z_i \rho = Z_i h + a_i \, .
\end{equation}
But the symplectic volume form, $\omega^d$, is, up to a constant
multiple, equal to
\begin{displaymath}
  e^h \prod \, dz_i \wedge \, d\bar{z}_i
  \prod \, dw_{\alpha} \wedge \, d\bar{w}_{\alpha} \, .
\end{displaymath}
Hence
\begin{displaymath}
  \div Z_i = L_{Z_i} \omega^d /\omega^d =Z_i h+1
\end{displaymath}
and by (\ref{eq:3.1})
\begin{displaymath}
  \phi_i = Z_i \rho \, ,
\end{displaymath}
so (\ref{eq:11.5}) reduces to
\begin{displaymath}
  \kappa \phi_{i} = \div Z_i + c_i \, .
\end{displaymath}
Note in particular that if the manifold, ``$M$'', in this section
is as in section~\ref{sec:7} the open dense subset, $M_{\stable}$,
of a \emph{compact} manifold, $M$, we can integrate
(\ref{eq:1.2}) over $M$, and using the fact that the integral
\begin{displaymath}
  \int (\div Z_i) \omega^d = \int \, d \iota (Z_i) \omega^d
\end{displaymath}
is zero, get
\begin{displaymath}
  \kappa \int \phi_i \omega^d = c_i \vol (M) \, .
\end{displaymath}
Thus if we normalize the moment map by requiring, as in
section~\ref{sec:7}, that the integral on the left be zero, we
end up with $c_i=0$, and the equation (\ref{eq:1.1}) reduces to
the simpler form
\begin{equation}
  \label{eq:11.6}
  \mu_{\lambda} - 2\sqrt{-1}  \partial \, \bar{\partial}
  \Log V_{\eff} =\kappa (\omega_{\lambda} + \sum \lambda_i
  \mu_i)\, .
\end{equation}

\noindent If we restrict to the case of level $0$, i.e., $\lambda =
0$, we reproduce an earlier result of Futaki's \cite{Fu}, Theorem
7.3.2, in the case of a K\"ahler--Einstein manifold $M$. Futaki treats
the case of $M$ a Fano manifold, which would be treated here in the
same fashion, {\em mutatis mutandis}.

\section{K\"ahler metrics on toric varieties}
\label{sec:12}

Let $\omega$ be the standard K\"ahler form on $\CC^n$:
\begin{displaymath}
  \omega = \sqrt{-1} \partial \, \bar{\partial} \rho
\end{displaymath}
where $\rho =\sum |z_i|^2$.  The linear action of $T=(S^1)^n$ on
$\CC^n$ is Hamiltonian with respect to this form, and by
(\ref{eq:3.1}), the moment map associated with the potential,
$\rho$, is
\begin{equation}
  \label{eq:12.1}
  \Phi (z) = (|z_1|^2 , \ldots , |z_n|^2) \, .
\end{equation}
Now let $G$ be a subtorous of $T$ and $\fg$ its Lie algebra.
{}From the inclusion of $G$ in $T$, one gets an inclusion of $\fg$
into $\RR^n$ and a transpose map
\begin{equation}
  \label{eq:12.2}
  L:\RR^n \to \fg^* \, ;
\end{equation}
and the moment map associated with the action of $G$ on $\CC^n$
is just the composition of (\ref{eq:12.1}) and (\ref{eq:12.2}).
It is easy to see that the images of the standard basis vectors,
$e_i$, of $\RR^n$ are the weights, $\alpha_i$, of the
representation of $G$ on $\CC^n$, so by (\ref{eq:12.1}) and
(\ref{eq:12.2}) this $G$-moment map can be written

\begin{equation}
\label{eq:12.3}
\psi (z) = \sum |z_i|^2 \alpha_i \, .
\end{equation}
This is a  proper mapping if and only if the $\alpha_i$'s are
polarized, i.e.,~if and only if there exists a $\xi \in \fg$
such that $\alpha_i (\xi)>0$ for all $i$.  We will assume
henceforth that this condition holds.

Now let $\lambda$ be a regular value of $\psi$ and assume that
$G$ acts freely on the level set, $\psi^{-1}(\lambda)$.  Then the
reduced space, $M_{\lambda} = \psi^{-1}(\lambda)$ is a K\"ahler
manifold with K\"ahler form, $\omega_{\lambda}$.  Moreover, since
the action of $T$ on $\CC^n$ commutes with the action of $G$, one
gets an induced Hamiltonian action of $T$ on $M_{\lambda}$ with
moment map, $\Phi_{\lambda} : M_{\lambda} \to \RR^n$.  Moreover,
if $p$ is the projection of $\psi^{-1} (\lambda)$ onto
$M_{\lambda}$
\begin{equation}
  \label{eq:12.4}
  \Phi_{\lambda} \circ p = \Phi \circ \iota_{\lambda}
\end{equation}
and this equation completely determines $\Phi_{\lambda}$.

The \emph{moment polytope} of the Hamiltonian $T$-manifold,
$M_{\lambda}$, is by definition the \emph{image} of $M_{\lambda}$
with respect to this moment map.  By (\ref{eq:12.4}) this is
equal to the image with respect to $\Phi$ of the set,
$\psi^{-1}(\lambda)$; and by (\ref{eq:12.1}) and (\ref{eq:12.2})
this is equal to the set
\begin{equation}
  \label{eq:12.5}
  \Delta_{\lambda} = \{ (t_1 , \ldots ,t_n) \in \RR^n_+ \, ,
    \sum t_i \alpha_i = \lambda \} \, .
\end{equation}
In other words it is the intersection
\begin{displaymath}
  \Delta_{\lambda} = \RR^n_+ \cap L^{-1} (\lambda)
\end{displaymath}
of the positive orthant in $\RR^n$ with the affine subspace,
$L^{-1}(\lambda)$.  It is easy to see that $\lambda$ is a regular
value of $\psi$ if and only if $L^{-1} (\lambda)$ intersects the
faces of $\RR^n_+$ transversally.  Thus the codimension $\kappa$
faces of $\Delta_{\lambda}$ are the intersections of
$L^{-1}(\lambda)$ with the codimension $\kappa$ faces of
$\RR^n_+$.  In particular the facets (codimension one faces) are
the intersections of $L^{-1}(\lambda)$ with the facets, $t_i=0$,
of $\RR^n_+$; so if $j:\Delta_{\lambda}\to\RR^n$ is the inclusion
map, the function
\begin{equation}
  \label{eq:12.6}
  j^* t_i =: \ell_i
\end{equation}
is just the distance function on $\Delta_{\lambda}$ to the
$i$\st{th} facet.  Now let $U$ be the open dense subset 
\begin{displaymath}
  (\CC^*)^n \cap \psi^{-1} (\lambda)/T
\end{displaymath}
of $M_{\lambda}$ and fix a base point, $c=(c_1, \ldots ,c_n)$, in
$\Delta_{\lambda}$.

\begin{theorem}
  \label{th:12.1}
The restriction of $\omega_{\lambda}$ to $U$ is equal to
$\sqrt{-1} \partial \, \bar{\partial} \rho_{\lambda}$ where
$\rho_{\lambda}$ is the pull-back by $\Phi_{\lambda}$ of the
function
\begin{equation}
  \label{eq:12.7}
  \sum \ell_i -c_i \Log \ell_i \, .
\end{equation}

\end{theorem}

\begin{proof}
  On $(\CC^*)^n$ the potential function 
  \begin{displaymath}
    \rho_1 = \sum |z_i |^2 -c_i \Log |z_i|^2
  \end{displaymath}
satisfies $\sqrt{-1}\partial \, \bar{\partial} \rho_1=\omega$;
and by (\ref{eq:3.3}), the $T$-moment map associated with this
potential is $\Phi -c$.  Thus the $G$-moment map associated with
this potential is $L$ composed with $\Phi -c$, or $\psi -
\lambda$.  In particular, the zero level set of this moment map
is the $\lambda$-level set of $\psi$ and so, by (\ref{eq:4.6}),
one gets the induced K\"ahler form on $M_{\lambda}$ by applying
$\sqrt{-1} \partial \, \bar{\partial}$ to the funciton
\begin{displaymath}
  \iota^*_{\lambda} (\sum |z_i|^2 - c_i \Log |z_i|^2)
\end{displaymath}
and by (\ref{eq:12.4}) and (\ref{eq:12.6}) this is just the
function
\begin{displaymath}
  \Phi^*_{\lambda} (\sum \ell_i - c_i \Log \ell_i) \, .
\end{displaymath}

\end{proof}

As we mentioned in the introduction, the proof above of
theorem~\ref{th:12.1} is a  slightly simplified version of the
proof described in \cite{CDG}.  For applications of this theorem,
see \cite{Ab} and \cite{CDG}.


\begin{thebibliography}{GGK2}

\bibitem[AM]{AM}
  R. Abraham and J. Marsden, Foundations of Mechanics,
  Benjamin--Cummings, Reading, MA (1978).

\bibitem[Ab]{Ab}
  Miguel Abreu, ``K\"ahler geometry of toric varieties and
  extremal metrics'', \emph{Int. J. Math.} \textbf{9} (1998),
  641--651.

\bibitem[At]{At}
  M. Atiyah, ``Convexity and commuting Hamiltonians'',
  \emph{Bull. Lond. Math. Soc.} \textbf{14} (1981), 1--15.

\bibitem[Be]{Be}
  A.L.Besse, Einstein Manifolds, \emph{Ergeb. Math. Grenzgb.}
  volume~\textbf{10}, Springer Verlag, Berlin, Heidelberg, New York
(1987). 


\bibitem[BG]{BG} 
  O. Biquard and P. Gauduchon, ``Hyper--K\"ahler metrics on
  cotangent bundles of Hermitian symmetric spaces'', \emph{Lecture notes
  in pure and applied mathematics} \textbf{184} Dekker, Amsterdam
(1997).

\bibitem[BD]{BD} D. Burns and P. DeBartolomeis, ``Stability of vector
bundles and extremal metrics'', \emph{Invent. Math.} \textbf{92}
(1988), 403-407.

\bibitem[CDG]{CDG}
  D. Calderbank, L. David and P. Gauduchon, ``About the Guillemin
  formula for the K\"ahler potential of a toric manifold'',
  preprint.

\bibitem[DH]{DH}
J.J. Duistermaat and G.L. Heckman, ``On the variation in
cohomology of the symplectic form of the reduced phase spce'',
\emph{Invent. Math.} \textbf{69} (1982), 259--269.

\bibitem[Fu]{Fu} A. Futaki, K\"ahler--Einstein Metrics and Integral
Invariants, \emph{Lecture Notes in Mathematics} \textbf{1314},
Springer Verlag, Berlin, Heidelberg, New York (1988). 

\bibitem[Gu]{Gu}
  V. Guillemin, ``K\"ahler structures on toric varieties'',
  \emph{J. Diff. Geom.} \textbf{40} (1984), 285--309.

\bibitem[HH]{HH} 
P. Heinzner and A. Huckleberry, ``K\"ahlerian
potentials and convexity properties of the moment map'',
\emph{Invent. Math.} \textbf{126} (1996), 65-84.

\bibitem[L]{L} 
C. LeBrun, ``Explicit self-dual metrics on
$\CC\mathbb{P}^2 \# \cdots \CC\mathbb{P}^2$'', \emph{J. Diff. Geom.}
\textbf{34}(1991), 223--253.

\bibitem[MFK]{MFK}
  D. Mumford, J. Fogarty and F. Kirwan, Geometric Invariant
  Theory, Springer, Berlin, Heidelberg, New York (1991).

\bibitem[PP]{PP}
H. Pedersen and Y.-S. Poon, ``Hamiltonian constructions of
K\"ahler-Einstein metrics and K\"hler metrics of constant scalar
curvature'', \emph{Comm. Math. Phys.} \textbf{136} (1991), 309--326.

\end{thebibliography}
\end{document}